\documentclass[12pt]{elsarticle}
\usepackage[usenames]{color}
\usepackage{amssymb}
\usepackage{amsmath}
\usepackage{amsthm}
\usepackage{amsfonts}
\usepackage{amscd}
\usepackage{graphicx}
\usepackage[colorlinks=true,
linkcolor=webgreen,
filecolor=webbrown,
citecolor=webgreen]{hyperref}
\definecolor{webgreen}{rgb}{0,.5,0}
\definecolor{webbrown}{rgb}{.6,0,0}
\usepackage{color}
\usepackage{fullpage}
\usepackage{float}
\usepackage{graphics}
\usepackage{latexsym}
\usepackage{epsf}
\usepackage{bookmark}
\usepackage{mathrsfs}
\usepackage{hyperref}
\def\C{\mathbb{C}}

\def\N{\mathbb{N}}

\def\cB{\mathcal B}
\def\cD{\mathcal D}
\def\cH{\mathcal H}
\def\le{\leqslant}
\def\ge{\geqslant}
\def\hc{\hat c}
\def\pr{^\prime }

\def\eop{\unskip\nobreak\hfil\penalty50\hskip2em\hbox{}\nobreak
\hfill\mbox{$\Box $}\par}
\newtheorem{theorem}{Theorem}[section]

\newtheorem{corollary}{Corollary}[section]
\begin{document}
\begin{frontmatter}
	\title{Discrete orthogonality of  hypergeometric polynomial sequences  on linear and quadratic lattices}
	
\author{Luis Verde-Star}
 \address{
Department of Mathematics, Universidad Aut\'onoma Metropolitana, Iztapalapa,
Apartado 55-534, Mexico City 09340,
 Mexico }
\ead{verde@xanum.uam.mx}

	\begin{abstract} 
	We present a  method to obtain   weight functions associated with linear and quadratic lattices that yield discrete orthogonality with respect to a quasi-definite moment functional  of  the orthogonal polynomial sequences in the Askey scheme, with the exception of the Jacobi, Bessel, Laguerre, and Hermite polynomials.


\end{abstract}
\end{frontmatter}

\section{Introduction}
The families of hypergeometric  orthogonal polynomial sequences included in the Askey scheme \cite{Hyp}  are some of the most important  and have been studied for a long time. Some of those families have an  orthogonality  determined by a weight function $w(x_k)$ defined on a finite or infinite  set of points $x_k$ that belong to a linear or quadratic lattice. In such cases the orthogonality of a polynomial sequence $\{u_n(t): n \in \N\}$ is given by 
\begin{equation}\label{discrOrtho}
\sum_{k \in  M} u_n(x_k) u_m(x_k) w(x_k)=  K_n \delta_{n,m}, \qquad  n,m \in M,
\end{equation}
where $M=\{0,1,2,\ldots,N\}$, for some positive integer $N$, or $M=\N$, and  $K_n$ is a  nonzero constant for $n\in M$.

The polynomial families that satisfy certain types of differential or difference equations and have a discrete orthogonality with positive weights $w(x_k)$ have been characterized and studied in great detail. The main references for discrete orthogonal polynomials are the books \cite{Niko} and \cite{Hyp}. See also \cite{AlvGMar}, \cite{DomMar1}, \cite{DomMar2}, and  \cite{SBFArea}. 

In the present paper we present a method to find a  weight function for each polynomial family in a subset $\cB$ of the Askey scheme that determines  a discrete orthogonality with respect to a quasi-definite moment functional. The families that are not included in $\cB$ are the ones of  Jacobi, Bessel, Laguerre, and Hermite polynomials. 
The weights $w(x_k)$ are the values at $t=1$ of a sequence of hypergeometric functions $f_k(t)$ that are of type $_3F_2$ or type $_2F_1$. The parameters of $f_k(t)$ depend on the parameters that determine the coefficients in the  three-term recurrence relation of the corresponding polynomial family. We show that the hypergeometric functions $f_k(t)$ are convergent at $t=1$ and that they are obtained from  $f_0(t)$ by repeated differentiation and multiplication by certain factors.
In this paper we do not deal with the problem of characterizing the cases for which the weights are positive.

For the classical families of discrete orthogonal polynomials of the Hahn, Meixner,  Krawtchouk, and Charlier polynomials, the orthogonality that we obtain coincides with the classical ones presented in \cite{Hyp}, with a different normalization in some cases.

We obtain our results using the linear algebraic approach of our previous papers \cite{Mops}, \cite{Rec}, \cite{PSGH}, and \cite{Uni}.
In  \cite{Uni} we presented a unified construction of all the  hypergeometric and basic hypergeometric  orthogonal polynomial sequences that uses three linearly recurrent sequences  of numbers that satisfy certain difference equation of order three, with constant coefficients. The initial terms of such numerical  sequences determine  the sequences of orthogonal polynomials and provide us with a uniform parametrization of all the  hypergeometric and basic hypergeometric sequences. In the present paper we used several results from \cite{Uni}.

In Section 2 we present some results related with  the construction of the hypergeometric orthogonal polynomials and their uniform parametrization from \cite{Uni}. In Section 3 we present the main result, whose proof uses infinite systems of linear equations and properties of hypergeometric functions. In Section 4 we find the weight functions for some families of orthogonal polynomials from the Askey scheme. Analogous results for the families of basic hypergeometric polynomial families in the $q$-Askey scheme will be presented elsewhere.

\section{The class $\cH_1$ of hypergeometric orthogonal polynomial sequences}

In this section we present some results about the class $\cH_1$ of the hypergeometric orthogonal polynomial sequences that were obtained in \cite{Uni}. 

Consider the homogeneous  difference equation
\begin{equation}\label{eq:diffeq}
 s_{k+3} = 3 ( s_{k+2} -s_{k+1}) + s_k, \qquad k\ge -1.
\end{equation}
The characteristic polynomial  of  this equation  has $1$ as a root of multiplicity 3, and therefore the general solution is a quadratic polynomial in $k$. 	Let $x_k$, $h_k$, and $e_k$ be solutions of \eqref{eq:diffeq}. Then   

\begin{equation}\label{eq:xhe}
 x_k= b_0 + b_1 k + b_2 k^2, \qquad h_k = a_0 + a_1 k + a_2 k^2,\qquad  
 e_k= d_0+ d_1 k + d_2 k^2. 
\end{equation}

The sequence $x_k$ determines the Newtonian basis $\{v_n(t): n\ge 0\}$ of the complex vector space of polynomials in $t$, defined by 
\begin{equation}\label{eq:Newton}
	v_n(t)= (t-x_0)(t-x_1)\cdots (t-x_{n-1}), \qquad n \ge 1,
\end{equation}
and $v_0(t)=1.$
We define the sequence $g_k$ by
\begin{equation} \label{eq:g}
g_k= x_{k-1} (h_k- h_0) + e_k, \qquad k \ge 1, 
\end{equation}
and $g_0=0$. Therefore we must have $e_0=d_0=0$. In addition, we  suppose that $h_k \ne h_j$ if $k \ne j$,  and $g_k \ne 0$ for $k \ge 1$. 

Let $\cD$ be the linear operator on the space of polynomials defined  by
\begin{equation}\label{eq:operD}
\cD v_k = h_k v_k + g_k v_{k-1}, \qquad k \ge 0.
\end{equation}
	Since $g_0=0 $ we see that $\cD t^n = h_n t^n +$ polynomial of lower degree. For $n \ge 0$ let $u_n$ be a monic polynomial of degree $n$ which is an eigenfunction of $\cD$ with eigenvalue $h_n$. That is 
\begin{equation}\label{eq:eigenEq}
\cD u_k = h_k u_k, \qquad k \ge 0.
\end{equation}

The operator $\cD$ is a generalized difference operator, which in concrete examples becomes a second order differential operator or a difference  operator on a linear or quadratic lattice. In \cite{Uni} we showed that
\begin{equation}\label{eq:u}
 u_n(t) = \sum_{k=0}^n c_{n,k} v_k(t), \qquad n \ge 0,
\end{equation}
 where the coefficients $ c_{n,k}$ are given by
 \begin{equation}\label{eq:cnk}
c_{n,k}=\prod_{j=k}^{n-1} \frac{g_{j+1}}{h_n -h_j}, \qquad 0 \le k \le n-1,
 \end{equation}
and $c_{n,n}=1$ for $n \ge 0$. This expression for $u_n(t)$ was also  obtained by Vinet and Zhedanov in \cite{VZ} using a different approach. The idea of representing orthogonal polynomials in terms of a Newtonian basis was introduced by Geronimus in \cite{Ger}. 

The matrix of coefficients $C=[c_{n,k}]$ is an infinite lower triangular matrix and the coefficients of $u_n$ appear in the $n$-th row of $C$. Since $c_{n,n}=1$ for $n \ge 0$, $C$ is invertible. 
	Let $C^{-1}=[\hc_{n,k}]$ and define the polynomials  
\begin{equation}\label{eq:Newthk}
w_{n,k}(t)= \prod_{j=k}^{n-1} (t-h_j), \qquad 0 \le k \le n.
\end{equation}
Using divided differences we obtain	
\begin{equation} \label{Cinverse}
	 \hc_{n,k}=\frac{ \prod_{j=k+1}^n g_j}{ w_{n+1,k}\pr (h_k) }=\prod_{j=k+1}^n \frac{g_j}{h_k - h_j}, \qquad 0 \le k \le n-1, 
\end{equation}
	 and $\hc_{n,n}=1$ for $n \ge 0.$

The entries in the $0$-th column of $C^{-1}$ are given by
\begin{equation}\label{eq:GenMom}
 \hc_{k,0} =\prod_{j=1}^k \frac{g_j}{h_0-h_j}, \qquad k \ge 1, 
\end{equation}
and $\hc_{0,0}=1$. We denote them by $m_k=\hc_{k,0}$ for $ k \ge 0.$  They satisfy $m_0=1$ and
\begin{equation} \label{moments}
	\sum_{k=0}^n c_{n,k} m_k = 0, \qquad n \ge 1.
\end{equation}
Note that the sequence $m_n$ satisfies a recurrence relation of order one.

In \cite{Uni} we also proved that the polynomial sequence $u_n(t)$ satisfies a three-term recurrence relation of the form
\begin{equation}\label{eq:3term}
	u_{n+1}(t)= (t-\beta_n) u_n(t) - \alpha_n u_{n-1}(t), \qquad n \ge 1, 
\end{equation}
where the coefficients are given by
\begin{equation}\label{eq:beta}
	\beta_n= x_n + \frac{g_{n+1}}{h_n - h_{n+1}} -\frac{g_n}{h_{n-1} - h_n}, 
\end{equation}
and 
\begin{equation}\label{eq:alpha}
 \alpha_n = \frac{g_n}{h_{n-1} -h_n} \left(\frac{g_{n-1}}{h_{n-2} - h_n} - \frac{g_n}{h_{n-1} - h_n} + \frac{g_{n+1}}{h_{n-1}-h_{n+1}} +x_n - x_{n-1} \right). 
\end{equation}

By Favard's theorem \cite[Thm. 4.4]{Chi}, \cite{MarAlv}, if all the $\alpha_n$ are positive and the $\beta_n$ are real then  $\{u_n\}$ is orthogonal with respect to a positive-definite moment functional, and if all the $\alpha_n$ are nonzero then $\{u_n\}$ is orthogonal with respect to a  quasi-definite moment functional.

Let $\cB$ be the class of all the families of hypergeometric orthogonal polynomial sequences whose recurrence coefficients are determined by a sequence of nodes $x_k=b_0 + b_1 k + b_2 k^2$, where at least one of $b_1$ or $b_2$ is nonzero. The families in the Askey scheme that are not included in $\cB$ are the families of Jacobi, Bessel, Laguerre, and Hermite. For these families $x_k=x_0$ for $k \ge 0$. 

\section{Generalized moments and discrete orthogonality}

Let us suppose that the parameters that determine the sequences $x_k, h_k, $ and $e_k$ are such that 
$ h_k \neq h_j $ if $k \ne j$, and the coefficients $\alpha_n$, given by \eqref{eq:alpha}, are nonzero for $n \ge 1$. Then, by Favard's theorem there exists a unique quasi-definite moment functional $\tau$ on the space of polynomials such that the polynomial sequence $\{u_n(t): n \ge 0\}$ is orthogonal with respect to $\tau$, that is,
\begin{equation}\label{orthogonality}
	\tau(u_n(t) u_m(t))= K_n \delta_{n,m}, \qquad K_n \ne 0, \ \ n,m \in \N.
	\end{equation}

  Since $u_0(t)=1$  we have
	\begin{equation}\label{eq:genmom}
	\tau(u_n(t)) =\sum_{k=0}^n c_{n,k} \tau(v_k(t) )=0, \qquad n \ge 1,
	\end{equation}
	and by the uniqueness of the inverse of $C$ we must have  $\tau(v_k(t))=\hc_{k,0}=m_k$  for $k \ge 0$.
This means that the numbers $m_k$ are the generalized moments of $\tau$ with respect to the Newtonian basis $\{v_k(t): k\ge 0\}$. Note that $\tau(1)=m_0=1.$

\begin{theorem}
	Suppose that at least one of $b_1$ and $b_2$ is nonzero,  at least one of $a_1$ and $a_2$ is nonzero, and $\alpha_n\ne 0$ for $n\ge 1$. Then there exists a unique weight function $w$, defined on the nodes $x_k$,  such that for every polynomial $p(t)$ we have 
\begin{equation}\label{eq:discrTau}
	\tau(p(t))=\sum_{k=0}^\infty p(x_k) w(x_k).
\end{equation}
\end{theorem}
 
{\it Proof:} From \eqref{eq:xhe} we see that the hypothesis implies that $\{x_k\}$ and $\{h_k\}$ are sequences of pairwise distinct numbers, and that the polynomial sequence $\{u_n(t): n\ge 0\}$ is orthogonal with respect to the quasi-definite  moment functional $\tau$.

We will write $r_j=w(x_j)$ for $j \ge 0$. Since  $\tau(v_k(t))=m_k$  for $k \ge 0$, the numbers $r_j$ that we want to find must satisfy 
\begin{equation}\label{eq:infsyst}
	m_k=\tau(v_k(t))= \sum_{j=0}^\infty  v_k(x_j)\, r_j , \qquad k \ge 0. 
\end{equation}
This is an infinite system of linear equations. Let us denote by $P$ the matrix of coefficients. Then  $P=[ v_k(x_j)]$, where $j$ is the index for rows, $k$ is the index for columns and $j$ and $k$ are non-negative integers. From the definition of the polynomials $v_k(t)$ in \eqref{eq:Newton} we can see that $P$ is an infinite  lower triangular matrix and its entries in the main diagonal are 
$$ v_k(x_k)=(x_k-x_0) (x_k-x_1) \cdots (x_k- x_{k-1}),$$
and they are nonzero. Therefore $P$ is invertible. 

Using basic properties of divided differences with respect to the nodes $x_j$ it is easy to show that $P^{-1}$ is the lower triangular  matrix whose $(j,k)$ entry equals $1/v_{j+1}^\prime(x_k)$.
Therefore, from the system \eqref{eq:infsyst} we obtain
\begin{equation}\label{eq:weights}
	r_k=\sum_{j=k}^\infty \frac{m_j}{v_{j+1}^\prime(x_k)}, \qquad k\ge 0.
\end{equation}
We will show  that if the parameters $a_1,a_2, b_1, b_2, d_1,d_2$  satisfy certain conditions, then all the series  in \eqref{eq:weights} are convergent (or terminating) hypergeometric series.

Let us define the power series
\begin{equation}\label{eq:fkt}
	f_k(t)= \sum_{j=k}^\infty \frac{m_j}{v_{j+1}^\prime(x_k)} t^j, \qquad k \ge 0.
\end{equation}

The denominators in the previous equation are 
\begin{equation}\label{eq:vprime}
	v_{j+1}^\prime(x_k)= \prod_{i=0, i\ne k}^j (x_k-x_i), \qquad k \ge  j \ge 0.
\end{equation}
By substitution of \eqref{eq:vprime} and \eqref{eq:GenMom} in \eqref{eq:fkt}, writing all the sequences using the formulas \eqref{eq:xhe} and \eqref{eq:g}, we find  an expression for $f_k(t)$ in terms of the parameters $ a_1,a_2, b_0,b_1,b_2,d_1,d_2$.  For example, the first 3 terms of $f_0(t)$ are
\begin{equation}\label{eq:f0t}
\begin{split}
	f_0(t)= 1+& \frac{a_1 b_0 + a_2 b_0 + d_1+d_2}{(a_1+a_2)(b_1+b_2)}t + \\
		    & \frac{(a_1 b_0 + a_2 b_0 + d_1+d_2)((b_0+b_1+b_2) (a_1+2 a_2)+d_1+2 d_2)}{ (a_1+a_2)(a_1+2a_2) (b_1+b_2)  (b_1+2b_2) } \frac{t^2}{2!} +\cdots 
  \end{split}
\end{equation}
Note that $a_0$ does not appear here because the sequence of eigenvalues $h_k$  enters only trough the  differences $h_k - h_j$. Therefore the value of $a_0$ is not relevant and we can suppose that  $a_0=0$.

The series in \eqref{eq:f0t} looks like a hypergeometric function. We can confirm that it is indeed a hypergeometric function through  the  changes of parameters that we define next. We will consider four cases, corresponding to whether $a_2$ and $b_2$ are zero or nonzero. 

\medskip
{\bf Case 1.} If   $a_2 \ne 0 $ and $b_2 \ne 0$ we  introduce the parameters $p,r,y_1,y_2$ as follows
\begin{equation}\label{eq:paramchange}
	\begin{split}
		a_1= & (r-1) a_2, \\
		b_1= & (p-1) b_2, \\
		d_1= & a_2 ( b_2 ( (y_1-1) (y_2-1) (p+r-y_1-y_2-2) + (r-1) (p-2) ) - (r-1) b_0),\\
		d_2= & a_2 ( b_2( (y_1+y_2) (p+r-y_1-y_2-1) + y_1 y_2 -r ( p-1)) -b_0).
	\end{split}
\end{equation}
The new  parameters $y_1$ and $y_2$  can be expressed in terms of the  original parameters by solving  for $y_1$ and $y_2$ the last two equations in \eqref{eq:paramchange}.

Substitution of  $a_1,b_1, d_1,d_2$ in $f_0(t)$   using \eqref{eq:paramchange}  gives us
\begin{equation}\label{eq:f0t3F2}
	f_0(t)= F(y_1,y_2,p+r-y_1-y_2-1; r,p;t),
\end{equation}
where $F$ is the hypergeometric function of type $_3F_2$ defined by
\begin{equation}\label{eq:3F2}
	F(y_1,y_2,y_3; r,p;t) = \sum_{k=0}^\infty \frac{(y_1)_k (y_2)_k (y_3)_k }{ (r)_k (p)_k} \frac{t^k}{k!},
\end{equation}
where $(a)_k$  denotes the shifted factorial,  defined  by $(a)_0=1$ and 
$$ (a)_k= a ( a+1) (a+2) \cdots(a+k-1), \qquad  a \in \C, \  k \ge 1. $$

For any values of the parameters $p,r,y_1$ and $y_2$ the parametric excess  (the sum of the lower parameters minus the sum of the upper parameters) of the hypergeometric function in \eqref{eq:f0t3F2} is  equal to one. Therefore it is convergent at $t=1$, since a hypergeometric function converges at $t=1$ if the real part of the parametric excess is positive. See \cite[Thm. 2.1.2, p.62]{AAR} or \cite[p.45 ]{Sla}.

The functions $f_k(t)$, for $k\ge 1$, can also be expressed in terms of hypergeometric functions using the parameters  $p,r,y_1$ and $y_2$. In those functions the factors of the form  $p+i$ in the denominators appear in a less regular way than in $f_0(t)$. This is so because they come from the derivatives $v_{j+1}^\prime(x_k)$ and this produces jumps and shifts  in the index $i$ of the factors $p+i$.

We can  verify by direct computations that 
\begin{equation}\label{eq:fmthyp1}
	\begin{split}
		f_k(t) = (-1)^k & \left( \frac{p+2k-1}{p+k-1}\right) \frac{t^k}{k!} \  \times \\
		                 &  \qquad D_t^k  F(y_1,y_2,r+p-y_1-y_2-1; r , p+ k;t),\quad k\ge 0,
	\end{split}
\end{equation}
where $D_t$ denotes differentiation with respect to $t$.

For $k \ge 0$ the parametric excess of $ F(y_1,y_2,r+p-y_1-y_2-1; r , p+ k;t)$  is $k+1$, and thus the function converges at $t=1$. Therefore $f_k(t)$ converges at $t=1$  and  the weights $r_k= f_k(1)$ are well defined for $k \ge 0.$ 

The functions $f_k(t)$ can also be expressed as  follows
\begin{equation}\label{eq:fmthyp2}
	\begin{split}
		f_k(t) = (-1)^k &  \frac{(y_1)_k (y_2)_k (p+r-y_1-y_2-1)_k}{(r)_k (p-1+k)_k} \frac{t^k}{k!} \  \times \\
		                 & \qquad  F(y_1+k,y_2+k,r+p-y_1-y_2-1+k; r+k , p+ 2 k;t).
	\end{split}
\end{equation}

The sequence of weights  $r_k=f_k(1)$  that we have found satisfies \eqref{eq:infsyst} for $k \ge  0$, and since the set  $\{v_k(t): k\ge 0\}$ is a basis for the space of all polynomials, we conclude that 
$\tau(p(t))=\sum_{k=0}^\infty  p(x_k) r_k$ for every polynomial $p(t)$.

\medskip
{\bf Case 2.} We consider now the case with $b_2=0$ and $a_2\ne 0$. We use here the substitutions
\begin{equation}\label{eq:subb2}
	\begin{split}
	b_2= & 0,\\
	a_1= & (r-1) a_2, \\
	d_1= & a_2 ( b_1 ( r-1+ (y_1-1) (y_2-1)) - b_0 (r-1)),\\
	d_2= & a_2 ( b_1 ( y_1+y_2 -r) -b_0), 
	\end{split}
\end{equation}
and obtain
\begin{equation}\label{eq:f0tb2}
	f_0(t)= F_2(y_1,y_2;r;t),
\end{equation}
where $F_2$ is the hypergeometric function of type $_2F_1$ defined by 
\begin{equation}\label{eq:F2}
	F_2(y_1,y_2;r;t)= \sum_{k=0}^\infty \frac{(y_1)_k (y_2)_k}{(r)_k} \frac{t^k}{k!}.
\end{equation}

From \eqref{eq:subb2} we can obtain $y_1$ and $y_2$ by solving a system of equations. They satisfy
$$y_1+y_2=\frac{(b_1 r + b_0) a_2 + d_2}{b_1 a_2},$$
and therefore the function \eqref{eq:f0tb2} converges at $t=1$ if the real part of  $ -(a_2 b_0 + d_2)/(a_2 b_1)$ is positive. 

The substitution of \eqref{eq:subb2} in $f_k(t)$ gives us 
\begin{equation}\label{eq:fkb2}
	f_k(t)= (-1)^k  \frac{t^k}{k!}\,  D_t^k F_2(y_1, y_2; r;t), \qquad k\ge 0.
\end{equation}

These functions are also  convergent at $t=1$ if the real part of $ -(a_2 b_0 + d_2)/(a_2 b_1)$ is positive. We can also write the functions $f_k(t)$ as
\begin{equation}\label{eq:fkb2bis}
	f_k(t)= (-1)^k \frac{(y_1)_k (y_2)_k}{(r)_k} \frac{t^k}{k!} F_2(y_1+k, y_2+k; r+k;t), \qquad k\ge 1.
\end{equation}

\medskip
{\bf Case 3.} Suppose now that  $a_2=0$  and $b_2 \ne 0$. In this case we see from \eqref{eq:f0t} that $f_0(t)$ becomes a simpler hypergeometric function. The parameter substitutions
\begin{equation}\label{eq:suba20}
	\begin{split}
		a_2= & 0,\\
		b_1= & (p-1) b_2, \\
		d_1= & a_1 ( b_2 ((y_1-1) (y_2-1) +p -2) -b_0),\\
		d_2= & -a_1 b_2 (p-y_1 -y_2 -1),
	\end{split}
\end{equation}
 in $f_0(t)$  give us 
\begin{equation}\label{eq:f0ta20}
	f_0(t)=F_2(y_1,y_2;p;t).
\end{equation}

The  parameters $y_1$ and $y_2$ can be obtained from \eqref{eq:suba20} by solving a system of equations. They  satisfy 
\begin{equation}\label{eq:sumy1y2}
	y_1+y_2= p - 1 +\frac{d_2}{a_1 b_2}.
\end{equation}
Therefore the parametric excess of the function   $F_2(y_1,y_2;p;t)$ in  \eqref{eq:f0ta20}   is 
$p- y_1-y_2 = 1-\frac{d_2}{a_1 b_2},  $
 and thus  $f_0(t)$ converges at $t=1$ if 
 \begin{equation}\label{eq:convF2}
	 \hbox{Re}\left( 1-\frac{d_2}{a_1 b_2}\right) > 0.
 \end{equation}

 In this case, applying the  substitutions of parameters \eqref{eq:suba20} to the series $f_k(t)$,  we obtain hypergeometric functions similar to   $F_2(y_1,y_2;p;t)$. By a straightforward computation we can verify that 
 \begin{equation}\label{eq:fkta20}
	 f_k(t)= (-1)^k \left( \frac{p+2k-1}{p+k-1}\right) \frac{t^k}{k!}\, D_t^k F_2(y_1,y_2;p+ k;t), \qquad k\ge 0.
 \end{equation}
 Note that all the functions $f_k(t)$ are convergent at $t=1$ if \eqref{eq:convF2} holds. Another expression for $f_k(t)$ in this case is the following
 \begin{equation}\label{eq:fkta20bis}
	 f_k(t)= (-1)^k  \frac{(y_1)_k (y_2)_k}{(p-1+k)_k} \frac{t^k}{k!} F_2(y_1+k,y_2+k;p+2 k;t).
 \end{equation}

\medskip
{\bf Case 4.} In this case we have $a_2=0$ and $b_2=0$ and  $f_0(t)$ has the form
 \begin{equation}\label{eq:f0a20b20}
	 f_0(t)=1+\frac{a_1 b_0 + d_1+d_2}{a_1 b_1}t + \frac{(a_1 b_0 + d_1+d_2)( a_1 b_0 +a_1 b_1 +d_1+2 d_2)}{a_1^2 b_1^2} \frac{ t^2}{2!} +\cdots
 \end{equation}

 The  substitutions 
 \begin{equation}\label{eq:suba2b2}
	 \begin{split}
		 a_2= & 0, \qquad b_2=0, \\
		 d_1= & a_1 (b_1(1+(y-1)z) -b_0),\\
		 d_2= & (z-1) a_1 b_1,
	 \end{split}
 \end{equation}
 give us 
 \begin{equation}\label{eq:f0a2b2}
	 f_0(t)=F_2(y,1;1;zt).
 \end{equation}

 By the ratio test we see that $F_2(y,1;1;z t)$ converges at $t=1$ if $|z|<1$.
With the change of parameters \eqref{eq:suba2b2} the series $f_k(t)$ become the hypergeometric functions
\begin{equation}\label{eq:fka2b2}
	f_k(t)= (-1)^k \frac{t^k}{k!}\, D_t^k(F_2(y,1;1; z t)), \qquad k\ge 1,
\end{equation}
where $D_t$ denotes differentiation with respect to $t$. It is clear that $f_k(t)$  converges at $t=1$ if  $|z|<1$.

From  \eqref{eq:suba2b2}  we can solve for $y$ and $z$ in terms of the other parameters and we obtain
\begin{equation}\label{eq:yza2b2}
		y=  \frac{a_1 b_0 + d_1 + d_2}{a_1 b_1 + d_2},\qquad \ \ 
		z=  1 + \frac{d_2}{a_1 b_1}.
\end{equation}
Therefore the functions  $f_k(t)$ are convergent  at $t=1$ if  $|1 + \frac{d_2}{a_1 b_1}|<1$. 

We have proved that in each of  the four cases the weights $r_k= f_k(1)$ are well defined if the parameters satisfy certain conditions. \eop

It is easy to verify that in  the four cases considered in the previous theorem we have
\begin{equation}\label{eq:sumcoefsfk}
	\sum_{k=0}^n \hbox{coefficient of } t^n \ \hbox{in} \ f_k(t)= 0, \qquad n\ge 1,
\end{equation}
and therefore
\begin{equation}\label{eq:sumfkt}
	\sum_{k=0}^\infty f_k(t) =1,
\end{equation}
for every $t$ in the region of convergence of $f_0(t)$, and in particular $\sum_{k=0}^\infty r_k=1$.  This property of the weights  also follows from $\tau(v_0)=\tau(1)=m_0=1.$

Let us note that in the four cases the hypergeometric function  $f_0(t)$ determines to a large extent the sequence of functions $f_k(t)$ and their region of convergence.

\begin{corollary}
	With the hypothesis of the previous theorem, if the parameters
	satisfy the conditions of the  corresponding case then   the polynomials
	$u_k(t)$ satisfy the discrete orthogonality
	\begin{equation}\label{eq:discreteortho}
		\tau( u_n(t) u_m(t)) =\sum_{k=0}^\infty u_n(x_k) u_m(x_k) r_k =K_n \delta_{n,m},
	\end{equation}
	where the constants $K_n=\alpha_1 \alpha_2 \cdots \alpha_n$ are nonzero for $n \ge 0.$
\end{corollary}

\section{ Some examples of the  weights $r_k$ }
In this section we find the weights $r_k$ that correspond to some families of orthogonal polynomial sequences from the Askey scheme \cite{Hyp}.

For each of the four cases we can find the coefficients of the three-term recurrence relation of the corresponding  polynomial family by applying to the  general formulas for   $\beta_n$ and $\alpha_n$, given in \eqref{eq:beta} and \eqref{eq:alpha}, the parameter substitutions that we  used in the proof of Theorem 3.1.

\medskip
{\bf Case 1.} In the first case the parameter substitution \eqref{eq:paramchange} applied to   $\alpha_n$  yields
\begin{equation}\label{eq:Wilson}
	\begin{split}
		\alpha_n=&\frac{b_2^2\, n  (n+y_1-1)(n+y_2-1) (n-p+y_1+y_2)}{(2n+r-3) (2n+r-2)^2(2n+r-1)} \ \times \\
 & (n-2+r) (n+r-y_1-1) (n+r -y_2-1) (n+p+r-y_1-y_2-2).
	\end{split}
\end{equation}
The change of parameters
\begin{equation}\label{eq:paramWils}
		r=  a+b+c+d,\quad p=2 c +1,\quad  y_1=c+d, \quad y_2=b+c, 
\end{equation}
transforms \eqref{eq:Wilson} into the sequence  $\alpha_n$ that corresponds to the Wilson polynomials written in terms of the parameters $a,b,c,d$ used in  \cite[eq. 9.1.5]{Hyp}.  The analogous result for the sequence $\beta_n$ of the Wilson polynomials holds. Therefore the weights $r_k$  for the Wilson polynomials are the values at $t=1$ of the functions $f_k(t)$ of \eqref{eq:fmthyp1} with the parameters changed by the substitutions \eqref{eq:paramWils}. We have
\begin{equation}\label{eq:fmWilson}
	\begin{split}
		f_k(t) =& (-1)^k  \left( \frac{2 c +2 k}{2 c +k}\right) \frac{t^k}{k!} \  \times \\
		  & \qquad D_t^k F(c+d,b+c,a+c;a+b+c+d, 2 c +1 +k; t).
	\end{split}
\end{equation}
For any values of $a,b,c,d$ the parametric excess of the hypergeometric function in the previous equation  is equal to $1+k$  and therefore  $f_k(t)$ converges at $t=1$ for $k\ge 0$. 

Note that \eqref{eq:paramWils} allows us to express  the parameters $a,b,c,d$  in terms of $ p,r,y_1,y_2$. 
From \eqref{eq:fmthyp2}  we obtain another expression for the functions  $f_k(t)$ of the Wilson polynomials.
\begin{equation}\label{eq:fmWilson2}
	\begin{split}
		f_k(t) =& (-1)^k   \frac{(c+d )_k (b+c )_k (a+c )_k}{(a+b+c+d)_k (2 c + k)_k} \frac{t^k}{k!} \  \times \\
		  &\qquad   F(c+d+k,b+c+k, a+c+k;a+b+c+d+k,  2 c +1+ 2 k;t).
	\end{split}
\end{equation}

If $y_1$ or $y_2$ is a negative integer then the hypergeometric  functions $f_k(t)$  defined in \eqref{eq:fmthyp1} are terminating. In that case the corresponding  set of   orthogonal polynomials is finite. This is the case for the Racah polynomial family.
  
\medskip
{\bf Case 2.} We present next examples that correspond to the  second case in the proof of Theorem 3.1, that is,  with   $b_2=0$  and $a_2$   nonzero. 

The substitution \eqref{eq:subb2} applied to the general formula for  $\alpha_n$ gives us
\begin{equation}\label{eq:alfC2}
	\alpha_n=-\frac{b_1^2  n (n+r-2) ( n+y_1-1) (n+y_2-1) ( n+r-y_1-1) (n+r-y_2-1)}{(2n+r-2)^2 (2n+r-3( ( 2n +r-1)}.
\end{equation}
The substitutions 
\begin{equation}\label{eq:subsContH}
	\begin{split}
		b_1= & i,\qquad  r= a+b+c+d,\\
		y_1=&a+c,\qquad  y_2=a+d,
	\end{split}
\end{equation}
transform \eqref{eq:alfC2} into
\begin{equation}\label{eq:alfContH}
	\begin{split}
		\alpha_n= & \frac{n (n+a+b+c+d-2)(n+a+c-1)}{(2n +a+b+c+d-2)^2(2n+a+b+c+d-1)} \times \\    &\qquad  \qquad \frac{ (n+a+d-1)(n+b+c-1)(n+b+d-1)}{(2n+a+b+c+d-3)} .
	\end{split}
\end{equation}
This $\alpha_n$ corresponds to the continuous Hahn polynomials \cite[Eq. 9.4.4]{Hyp}.

Applying the change of parameters  \eqref{eq:subsContH} to the functions $f_k(t)$ in \eqref{eq:fkb2} we obtain
\begin{equation}\label{eq:rmContH}
	f_k(t)= (-1)^k  \frac{t^k}{k!}\, D_t^k F_2(a+c, a+d; a+b+c+d;t), \qquad k\ge 0.
\end{equation}
These hypergeometric functions converge at $t=1$ if the real part of $ b-a$ is positive.

\medskip
The Hahn polynomials \cite[Eq.9.5.4]{Hyp}, with parameters $\alpha, \beta, N$,  are obtained from the general formulas for $\alpha_n$ and $\beta_n$ with the substitutions
\begin{equation}\label{eq:subsHahn}
	\begin{split}
		b_2= & 0,\qquad b_1=1, \qquad b_0=0,\\
		a_1= & (\alpha+\beta+1) a_2, \\
		d_1= & -a_2 (N \alpha -\beta-1),\\
		d_2= & -a_2 (N +\beta +1).
	\end{split}
\end{equation}
With these values of the parameters we obtain $f_0(t)=F_2(-N, \alpha+1;\alpha+\beta+2;t)$, and this terminating  hypergeometric function determines the discrete orthogonality \cite[Eq. 9.5.2]{Hyp} of the Hahn polynomials.

\medskip
{\bf Case 3.} We consider next the  case with $a_2=0$ and $b_2\ne 0$.
The continuous dual Hahn polynomials, with parameters $a,b,c$ \cite[Eq.9.3,5]{Hyp}, are obtained from the general formulas taking
\begin{equation}\label{eq:contDualH}
	\begin{split}
	a_2=& 0,\qquad a_1=1, \\
	b_2=& 1, \qquad    b_1=2 c, \qquad  b_0=c^2,\\
   d_2= &a+b, \qquad  	d_1= a b+ a c + b c -a-b.
	\end{split}
\end{equation}
Combining these equations with \eqref{eq:suba20} we obtain
\begin{equation}\label{eq:contDualHpy}
	p= 2 c +1,\qquad y_1=a+c, \qquad y_2=b+c.
\end{equation}

By substitution in  \eqref{eq:fkta20} we obtain
 \begin{equation}\label{eq:fkContDualH}
	 f_k(t)= (-1)^k \left( \frac{2c+2k}{2 c +k}\right) \frac{t^k}{k!}\, D_t^k F_2(a+c,b+c;2c+1+ k;t), \qquad k\ge 0.
 \end{equation}
These functions are convergent at $t=1$ if the real part of $1-a-b$ is positive.
The numbers $r_k=f_k(1)$ are the weights for the continuous dual Hahn polynomials.

\medskip
{\bf Case 4.} We consider here examples with $a_2=0$ and $b_2=0$. Some of the classical discrete orthogonal polynomials fall in this case.

If  $a_2=0$ and $b_2=0$   the functions $f_k(t)$ are given in \eqref{eq:fka2b2}, which is
\begin{equation}\label{eq:fka2b2bis}
	f_k(t)= (-1)^k \frac{t^k}{k!}\, D_t^k\, F_2(y,1;1; z t), \qquad k\ge 0,
\end{equation}
where $y$ and $z$ are given in \eqref{eq:yza2b2}, and  all the $f_k(t)$ are convergent at $t=1$ if $|z|<1$. It is clear from \eqref{eq:fka2b2bis} that  $f_0(t)$ determines all the functions $f_k(t)$, and the weights $r_k=f_k(1)$.
	
The Krawtchouk polynomials \cite[Eq. 9.11.4]{Hyp}, with parameters $p$ and $N$, are obtained with
\begin{equation}
	\begin{split}
		a_2=&0, \qquad b_2=0, \qquad b_1=1,\\
		d_2=& (p-1) a_1, \qquad d_1=-a_1 ( Np + p +b_0 -1),\\
	\end{split}
\end{equation} 
and these equations give $f_0(t)=F_2(-N,1;1;pt)$, which yields the discrete orthogonality in \cite[Eq. 9.11.2]{Hyp}.

The Meixner polynomials, with parameters $c$ and $\beta$  \cite[Eq. 9.10.4]{Hyp}, are obtained with
\begin{equation}
	\begin{split}
		a_2=&0, \qquad b_2=0, \qquad b_1=1,\\
		d_2=&\frac{a_1}{c-1}, \qquad d_1=\frac{a_1(\beta -b_0) c +b_0-1)}{c-1},\\
	\end{split}
\end{equation} 
and these parameters give us  $f_0(t)=F_2(\beta,1;1;\frac{ct}{c-1})$, which produces the discrete orthogonality described in \cite[Eq. 9.10.2]{Hyp}.

The Charlier polynomials, with parameter $a$ \cite[Eq. 9.14.4]{Hyp}, are obtained with
\begin{equation}
	\begin{split}
		a_2=&0, \qquad b_2=0, \qquad b_1=1,\\
		d_2=& -a_1, \qquad d_1=  (1-a -b_0) a_1,\\
	\end{split}
\end{equation} 
which give $f_0(t)= \exp(-a t).$ This function produces the discrete orthogonality \cite[Eq. 9.14.2]{Hyp}.

\section{Conclusions.}
We have shown that all  the hypergeometric orthogonal polynomial sequences in the Askey scheme, other than the very classical families (Jacobi, Bessel, Laguerre, and Hermite),  satisfy a discrete orthogonality on a sequence of distinct  nodes $x_k$ in the complex plane. 
The nodes are of the form $x_k=b_0 + b_1 k + b_2 k^2$, where at least one of $b_1$ and $b_2$ is nonzero. The weights $r_k$ associated with the nodes $x_k$ are the values at $t=1$ of a sequence of hypergeometric functions. We have not dealt with the problem of characterizing the cases for which the weights $r_k$ are positive. 

The approach used in this paper can be used to obtain   discrete orthogonality  for some of the families of basic hypergeometric  orthogonal polynomials in the $q$-Askey scheme.

\end{document}